\documentclass[12pt]{article}
\usepackage{amssymb}
\makeatletter
\date{}

\newtheorem{theorem}{Theorem}[section]

\newtheorem{proposition}[theorem]{Proposition}

\newcommand{\z}{{\Bbb Z}}
\newcommand{\q}{{\Bbb Q}}
\newcommand{\re}{{\Bbb R}}
\newcommand{\N}{{\Bbb N}}

\newcommand{\invlim}{{\rm invlim}}

\newcommand{\lo}{\longrightarrow}

\newcommand{\black}{{\blacksquare}}
\newcommand{\tor}{{\rm Tor}}
\newcommand{\diam}{{\rm diam}}
\newcommand{\dist}{{\rm dist}}

\begin{document}

\title{On dimensionally exotic maps}

\author{ Alexander Dranishnikov\footnote{ the first author was supported by NSF grant  DMS-0904278;}  \ and Michael Levin\footnote{the second author was supported by ISF grant 836/08}}

\maketitle
\begin{abstract} 
We call a value $y=f(x)$ of a map $f:X\to Y$ dimensionally regular
if $\dim X\le \dim(Y\times f^{-1}(y))$.
It was shown in \cite{first-exotic} that if a map $f:X\to Y$ between compact metric spaces
does not have dimensionally regular values, then $X$ is a Boltyanskii compactum, i.e. a compactum satisfying the
equality $\dim(X\times X)=2\dim X-1$.
In this paper we prove that every Boltyanskii compactum $X$ of dimension $\dim X \geq 6$ admits a map
$f:X\to Y$ without dimensionally regular values. Also we exhibit a 4-dimensional Boltyanskii compactum
for which every map has a dimensionally regular value.

\

\bigskip
\\
{\bf Math. Subj. Class.:}  55M10 (54F45 55N45)
\end{abstract}
\begin{section}{Introduction}\label{intro}
Throughout this paper we assume that  maps are continuous and
spaces are separable metrizable. We recall that a {\em compactum}
means a compact metric space. By dimension $\dim X$ of a space $X$  we
assume the covering dimension.

The  famous Sard theorem states that a smooth map of closed manifolds $f:X\to Y$ has a regular value.
How far Sard's theorem can be extended for mappings of compact metric spaces? We study this question from
dimension theoretic point of view. We call a value $y\in Y$  for a continuous mapping of compact metric spaces
$f:X\to Y$ to be {\em dimensionally regular} if $$\dim X\le \dim(Y\times f^{-1}(y)).$$ Note that a regular value of a smooth map
is dimensionally regular. The question we study is whether every continuous map between compacta has a dimensionally regular value. 

In our previous paper~\cite{first-exotic} we proved that this question has an affirmative answer if compactum $X$ has the property $\dim(X\times X)=2\dim X$. It is known that not all compact metric spaces are such. 
The first example of a compactum which do not satisfies this equality was constructed by Boltyanskii~\cite{Bolt}. His example was 2-dimensional with the dimension of the square equal 3. We call  a compactum $X$  a {\em Boltyanskii compactum} if it has the property $\dim (X\times X)\ne 2\dim X$.
It is known that for Boltyanskii compacta necessarily $\dim(X\times X)= 2\dim X-1$.

One of the main discovery of~\cite{first-exotic} is the existence of continuous maps $f:X\to Y$ without dimensionally
regular values. We called such maps {\em dimensionally exotic}. We proved the following
\begin{theorem}
{\rm (\cite{first-exotic})}
\label{particular-exotic-map}
For every $n \geq 4$ there is an $n$-dimensional Boltyanskii compactum $X$
admitting a dimensionally exotic map $f : X \lo Y$ to a $2$-dimensional compactum $Y$.
\end{theorem}

The question remained unanswered whether every Boltyanskii compactum admits an exotic map.
We address it in this paper. Namely, we prove the following:
\begin{theorem}
\label{general-exotic}
Every  finite dimensional Boltyanskii compactum $X$ with $\dim X \geq 6$ admits
a dimensionally exotic map $f : X \lo Y$ to a $4$-dimensional compactum $Y$.
\end{theorem}

In~\cite{first-exotic} we observed that every continuous map of a compactum of dimension $\le 3$ has 
dimensionally regular values. It turns out that in the case of dimension 4  there is no  uniform answer.
\begin{theorem}
\label{dim=4}
There is a $4$-dimensional Boltyanskii compactum $X$ such that every continuous map $f:X\to Y$ has a dimensionally
regular value.
\end{theorem}
 It still remains open whether every Boltyanskii compactum of dimension  $5$ admits a dimensionally exotic map.
We note that the proof of Theorem~\ref{general-exotic} closely related to the construction by the first author
of 4-dimensional ANR compacta $X$ and $Y$ with  $\dim(X\times Y)=7$~\cite{DrUsp},\cite{DrArxiv}. It is a long
standing open problem whether such phenomenon can happen for 3-dimensional ANRs. Thus, our remaining problem
in dimension 5 could be closely related to this one and quite difficult.

\

  In the next section we briefly review
 basic facts of 
Cohomological Dimension   and prove  Theorem \ref{dim=4}.
Auxiliary  constructions and propositions  for proving 
Theorem \ref{general-exotic} are presented in Section \ref{auxiliary}
and Theorem  \ref{general-exotic} is proved in the last section.

 \end{section}

\begin{section}{Cohomological Dimension }
\label{cohomological-dimension} In this section we review basic facts of Cohomology Dimension
and prove Theorem \ref{dim=4}.

By   cohomology we always mean the Cech
cohomology. Let $G$ be an abelian group. The {\bf cohomological dimension} $\dim_GX$
of a space $X$ with respect to the coefficient group $G$ does not exceed $n$, $\dim _G X \leq n$ if $H^{n+1}(X,A;G)=0$ for every closed $A
\subset X$. We note that this condition implies that
$H^{n+k}(X,A;G)=0$ for all $k\ge 1$ \cite{Ku},\cite{DrArxiv}.
Thus, $\dim _G X =$ the smallest integer $n\geq 0$ satisfying
$\dim _G X \leq n$ and $\dim _G X = \infty $ if such an integer
does not exist. Clearly, $\dim_G X \leq \dim_{\z}X\le\dim X$.

 \begin{theorem}
 \label{alexandrov}
 {\rm ({\bf Alexandroff})}
$\dim X=\dim_\z X$ if $X$ is a finite dimensional space.
\end{theorem}

Let $\mathcal P$ denote the set of all primes. The {\em  Bockstein basis} is the collection of groups
$\sigma=\{ \q, \z_p , \z_{p^\infty}, \z_{(p)} \mid p\in\mathcal P
\}$ where $\z_p =\z/p\z$ is the $p$-cyclic group,
$\z_{p^\infty}={\rm dirlim} \z_{p^k}$  is the $p$-adic circle, and
$\z_{(p)}=\{ m/n \mid n$ is not divisible by $p \}\subset\q$  is the $p$-localization of integers.
\\\\
 The Bockstein basis   of an abelian group $G$ is the collection
$\sigma(G) \subset \sigma$ determined by the rule:

$\z_{(p)} \in \sigma(G)$ if $G/\tor G$ is not divisible by $p$;

$\z_p \in \sigma(G)$ if $p$-$\tor G$ is not divisible by $p$;

$\z_{p^\infty} \in \sigma(G)$ if $p$-$\tor G\neq 0$ is  divisible by $p$;

$\q \in \sigma(G)$ if $G/\tor G\neq 0$ is  divisible by all $p$.
\\
\\
Thus  $\sigma(\z)=\{\z_{(p)}\mid p\in\mathcal P\}$.

\begin{theorem}
\label{bockstein-theorem}
{\rm ({\bf  Bockstein Theorem})} For a compactum $X$,
$$\dim_G X =\sup \{ \dim_H X : H \in \sigma(G)\}.$$
\end{theorem}
${}$\\
Suggested by the Bockstein inequalities we say that
a function $D:\sigma\lo \N \cup \{0, \infty\}$ is   a {\em
$p$-regular} if
$
D(\z_{(p)})=D(\z_p)=D(\z_{p^{\infty}})=D(\q)$ and it  is   {\em
$p$-singular} if
$ D(\z_{(p)})=\max\{ D(\q),D(\z_{p^{\infty}})+1\}.$
A $p$-singular function $D$ is called $p^+$-singular if
$D(\z_{p^{\infty}})=D(\z_p)$ and it is
called {\em $p^-$-singular} if
$D(\z_{p^{\infty}})=D(\z_p)-1$. A function
$D:\sigma\lo \N \cup \{0, \infty\}$ is called a {\em dimension
type} if for every prime $p$ it is either $p$-regular or
$p^{\pm}$-singular. 
Thus, the values of $D(F)$ for the Bockstein fields $F\in\{\z_p,\q\}$ 
together with $p$-singularity types of $D$
determine the value $D(G)$ for all groups in $\sigma$.
For a dimension type $D$ denote $\dim D =\sup \{ D(G) : G \in \sigma \}$.

\begin{theorem}
\label{bockstein-inequalities}
{\rm ({\bf Bockstein Inequalities} \cite{Ku}, \cite{DrArxiv})}
For every space $X$ the function $d_X:\sigma\lo \N \cup \{0,
\infty\}$ defined as $d_X(G)=\dim_GX$ is a dimension type. 
\end{theorem}
If $X$ is compactum
$d_X$ is called the {\em dimension type of $X$}.

\begin{theorem}
\label{dranishnikov-realization-theorem}
{\rm ({\bf Dranishnikov  Realization Theorem} \cite{DrUsp},\cite{DrPacific})}
For every   dimension type $D$ there is a compactum $X$
with $d_X=D$ and $\dim X=\dim D$.
\end{theorem}

\begin{theorem}
{\rm({\bf Olszewski Completion Theorem}\cite{Ol})}
\label{completion}
For every space $X$ there is a complete space $X'$ such that
$X \subset X'$ and $d_X =d_{X'}$.
\end{theorem}

Let $D$ be a dimension type. We will use the abbreviations
 $D(0)=D(\q)$, $D(p)=D(\z_p)$. Additionally, if $D(p)=n\in\N$ we will write $D(p)=n^+$ if $D$ is $p^+$-regular
and $D(p)=n^-$ if it is $p^-$-regular. For a $p$-regular $D$  we leave it without decoration: $D(p)=n$.
Thus, any sequence of decorated numbers $D(p)\in\N$, where $p\in\mathcal P\cup\{0\}$  define a unique dimension type.
There is a natural order on decorated numbers $$\dots<n^-<n<n^+<(n+1)^-< \dots\ .$$
Note that the inequality of dimension types $D\le D'$ as functions on $\sigma$ is equivalent to the family of inequalities
$D(p)\le D'(p)$ for the above order for all $p\in\mathcal P\cup\{0\}$.
The natural involution on decorated numbers that exchange the decorations '+' and '-' keeping the base fixed defines
an involution $\ast $ on the set of dimension types . Thus, $\ast$ takes
$p^+$-singular function $D$ to $p^-$-singular $D^*$ and vise versa.

Let $D_1$ and $D_2$ are dimension types. 
Suggested by Bockstein Product Theorem we define the dimension type $D_1 \boxplus D_2$ as follows.
If $D_1(p)=n^{\epsilon_1}$ and $D_2(p)=m^{\epsilon_2}$
where $\epsilon_i$ is a decoration, i.e., '+' or '-' or empty, then 
 $$(D_1\boxplus D_2)(p)=(n+m)^{\epsilon_1\otimes\epsilon_2}$$
with the product of the signs $\epsilon_1\otimes\epsilon_2$ defined by
$$\epsilon\otimes empty=\epsilon,\ \  \ \ \epsilon\otimes \epsilon=\epsilon, \ \ \epsilon=\pm,\ \ \  and\ \ \
+\otimes -=-.$$
\\
It turns out that $D_1 \boxplus D_2$ is indeed a dimension type and its definition is justified by

\begin{theorem}
\label{bockstein-product-theorem}
{\rm(\bf {Bockstein Product Theorem} \cite{DrArxiv},\cite{Sch}, \cite{DyA})}
For  any two compacta  $X$ and $Y$ 

$$d_{X \times Y} =d_X \boxplus d_Y.$$
\end{theorem}

By $B_n$ we denote the dimension type such that
$B_n(p)=(n-1)^+$ for all $p\in\mathcal P$ and $B_n(\q)=n-1$.
The product formula implies that an $n$-dimensional compactum $X$ is a Boltyanskii compactum if and only if
$d_X \leq B_n$.

 For an integer  $n\geq 0$ we denote by $n$ the dimension type which sends every
 $G\in \sigma$ to $n$. For a dimension type $D$  by $D +n$ we mean  the ordinary sum
 of $D$ and $n$ as functions. Note that $D +n$ is also a dimension type.
 
 For dimension types $D_1$ and $D_2$
  we introduce another  important operation
 $$D_1 \oplus D_2=(D_1^*\boxplus D_2^*)^*. $$
  It can be shown   that $D_1 \oplus D_2$
 is a dimension type as well. 
  The operation $\oplus$  is justified by   the following properties.

 \begin{theorem}
 \label {union}
 {\rm ({\bf Dydak Union Theorem}\cite{Dy},\cite{first-exotic})}
 Let $X$ be a compactum and $D_1$ and $D_2$  dimension types and let 
 $X=A \cup B$ be  a decomposition with $d_A \leq D_1$ and
 $d_B \leq D_2$.  Then $d_X\leq D_1 \oplus D_2 +1$.
\end{theorem}

\begin{theorem}
\label{decomposition}
{\rm ({\bf Dranishnikov  Decomposition  Theorem} \cite{DrPacific}, \cite{first-exotic})}
  Let $X$ be a finite dimensional compactum 
  and $D_1$ and $D_2$ dimension types such that   $d_X \leq D_1 \oplus D_2 +1$.
 Then there is a decomposition  $X=A \cup B$ such that
 $d_A \leq D_1$ and   $d_B \leq D_2$.
\end{theorem} 

\begin{theorem}
\label{levin-lewis}
{\rm ({\bf Levin-Lewis} \cite{levin-lewis}, \cite{first-exotic})}
Let $f: X \lo Y$ be a map of compacta  such that
$d_f \leq D_1$ and $d_Y \leq D_2$ where 
$d_f \leq D_1$ means that the dimension type of each fiber of $f$
is less or equal $D_1$.
Then $d_X \leq D_1 \oplus D_2$.
\end{theorem}

Let $D_1$ and $D_2$ be dimension types.
It can be shown that $D_1 \boxplus D_2 \leq D_1 \oplus D_2$ 
 and if $D_1'$ are $ D_2'$ are dimension types such that $D_1 \leq D_1'$ and
$D_2 \leq D_2'$ then $D_1 \boxplus D_2 \leq D_1'\boxplus D_2'$ and
$D_1 \oplus D_2 \leq D_1' \oplus D_2'$.
 Note that $D \boxplus n =D \oplus n =D+n$ for any  dimension type $D$. 
 
 An $n$-dimensional space $X$
    is said to be dimensionally full-valued if $d_X=n$. 
        Note that every $0$ or $1$-dimensional
    compactum is dimensionally full-valued. P. S. Alexandroff proved the following dual statement:
    
    \begin{theorem}
    \label{co-dim=1}
    {\rm (\cite{Ku})}
    Let $X$ be an  $(n-1)$-dimensional compact
    subset of $\re^n$. Then $X$ is dimensionally  full-valued.
    \end{theorem}
    
    Now we are ready to prove
    \begin{theorem}
    \label{cube}
  Let $X$ be an $n$-dimensional Boltyanskii compactum, $n \geq 4$, containing 
  an $(n-1)$-dimensional cube. Then $X$ does not admit a dimensionally exotic map
  $f : X \lo Y$ to a   compactum $Y$ with $\dim Y\leq 2$.
  \end{theorem}
  {\bf Proof.} Let $B \subset X$ be an $(n-1)$-dimensional cube. Aiming
  at a contradiction assume that there is a dimensionally exotic map $f : X \to Y$
  with $\dim Y \leq 2$.  First notice that $Y$ is not dimensionally full-valued since
  otherwise $\dim (Y \times f^{-1}(y))=\dim Y +\dim f^{-1}(y)$ and 
  by Hurewicz's theorem $f$ cannot be exotic. Thus we have $\dim Y=2$.
  Denote  $k=\dim f|_B$.
  Again by Hurewicz's theorem $k \geq n-3$. If $k=n-3$ then by 
   Theorem \ref{levin-lewis} we have 
    $n-1 =d_B \leq d_Y \oplus (n-3)=d_Y + (n-3)$
    and this contradicts to the fact that $Y$ is $2$-dimensional and not
dimensionally full-valued. The case $k=n-1$ is also impossible because
then  there is a fiber $F$ of $f$ containing an $(n-1)$-dimensional cube
and hence  $\dim (Y \times F) =2 + n-1=n+1$  that contradicts the fact that  $f$ is exotic.
So we are left with the only case   $k=n-2$. Then there is a fiber $F$
of $f$ with $\dim F\cap B =n-2$ and, by Proposition \ref{co-dim=1},
$d_F \geq n-2$. Thus $d_{ Y \times F} =d_Y \boxplus d_F \geq d_Y +(n-2)$
and hence $\dim (Y\times F) \geq \dim Y +n-2=n$ and this contradicts 
the assumption that $f$ is dimensionally exotic. The theorem is proved.
$\black$. 
  \\\\
 {\bf Proof of  Theorem \ref{dim=4}}. It   follows immediately from Theorem \ref{cube}. 
 Indeed, take a $4$-dimensional
  Boltyanskii compactum $X$  containing a $3$-dimensional cube and assume that 
  $f : X \lo Y$ is a dimensionally exotic map to a compactum $Y$. Clearly 
  $\dim Y \leq 3$ and  by Theorem \ref{cube} we have
  $\dim Y = 3$. Then by Hurewicz's theorem there is a fiber $F$ of $f$
  with $\dim F \geq 1$ and hence $f$ cannot be exotic because 
  $\dim (Y \times F) \geq \dim Y +1 =4$. $\black$
  \\\\
  Note that Theorem \ref{cube} also implies that the dimension of $Y$
  in Theorem \ref{general-exotic} cannot be reduced to $2$. We don't know
  if it can be reduced to $3$.

\end{section}

\begin{section}{Auxiliary propositions and constructions}
\label{auxiliary}
Let $X \lo Y$ be a map to a simplicial complex $Y$ and $A \subset X$. We say
that $f$ is {\em dimensionally deficient} on $A$
if $\dim(f^{-1}(\Delta) \cap A) < \dim \Delta$ for every
simplex $\Delta$ in $Y$.
The following two propositions are straightforward applications of the Baire category theorem and their proof is left to the reader.

\begin{proposition}
\label{extension}
Let $X$ be a compactum, $M$ a triangulated  $n$-dimensional
manifold possibly with boundary,
$A$ a $\sigma$-compact subset of  $X$ such that $\dim A \leq n-1$,
$F$ a closed subset of $X$, and let $f : X \lo M$  be a map
which  is dimensionally
deficient on $A \cap F$. Then 
the map $f$ can be arbitrarily closely approximated by a map $f' : X \lo M$
such that $f'$ is dimensionally deficient on $A$ and $f'$ coincides with $f$ on $F$.
\end{proposition}

\begin{proposition}
\label{subdivision}
Let $A \subset X$ be a $\sigma$-compact subset of a compactum  $X$, and let
$f : X \lo K$  be a map  to a finite simplicial complex $K$
such that $f$ is  dimensionally deficient on $A$. Then
for every subdivision $K'$ of $K$ the map $f$ can be arbitrarily approximated
by a map $f' : X \lo K'$ such that $f'$ is dimensionally deficient  on $A$ with
respect to the triangulation of $K'$ and for every simplex $\Delta$ of $K$
we have that $f^{-1}(\Delta)=f'^{-1}(\Delta)$.
\end{proposition}

The following proposition besides the Baire category theorem  uses the fact $n+1$
$(n-1)$-dimensional planes in $\re^n$ in a general position do not intersect.
This fact  implies that a generic map $f:Y\to\Delta^n$ of an $(n-1)$-dimensional compactum
to an $n$-dimensional simplex has $|f^{-1}(x)|\le n$ for all $x\in\Delta^n$.

\begin{proposition}
\label{finite-to-one-approximation}
Let $X$ be a compactum, $A$ a $\sigma$-compact subset, and let $f: X \lo K$ be
a map to a finite simplicial complex $K$ such that $f$ is dimensionally deficient
on $A$. Then $f$ can be arbitrarily closely approximated by a map
$f'$ such that for every simplex $\Delta$ of $K$ and for every $y \in \Delta$
we have that $f^{-1}(\Delta)=f'^{-1}(\Delta)$ and the number of points in
 $f'^{-1}(y)\cap A$ does not exceed  $\dim \Delta$.
Clearly that $f'$ is dimensionally deficient  on $A$ as well.
\end{proposition}

 We will also need
\begin{proposition}
{\rm (\cite{first-exotic})}
\label{rational-dimension}
Let $X$ be a finite dimensional compactum and $n>0$.
Then  $\dim_\q  X \leq n$ if and only if for every closed subset
$A$ of $X$ and every map $f : A \lo S^n$ there is a map
$g : S^n \lo S^n$ of non-zero degree such that
$g \circ f : X \lo S^n$ continuously extends over $X$.
\end{proposition}

 Let $g,h : X \lo Y$ be maps. Consider the product $X \times [0,1]$
 and consider $g$ and $h$ as the maps from $X\times \{0\}$ and $X\times \{1\}$
 respectively.
 We recall that the {\em double mapping cylinder}  $M(g,h)$ of  $g$ and $h$
 is the quotient space of $X \times [0,1]$ in which
 $X\times \{0\}$ is identified with $Y$ according to the map $g$  and
 $X\times \{1\}$
 is identified with another copy of $Y$ according to the map $h$.
   The quotient map $P : X\times [0,1] \lo M(g,h)$ will be called
 the {\em cylinder projection} and
 the map $p: M(g,h) \lo [0,1] $ induced  by the map $(x,t)\lo t$
 will be called the {\em interval projection}. We call $p^{-1}(0)$ and $p^{-1}(1)$
 the $g$-part and the $h$-part of $M(g, h)$ respectively.
 For pointed maps $g,h : (X,x_0)\lo (Y,y_0)$  we identify the  interval  $[0,1]$
 with the interval $P(\{x_0 \} \times [0,1])$ and call it the {\em axis} of $M(g,h)$.
  We can  rescale or shift the interval $[0,1]$
 to an interval $[a,b]$ and in that case we say that $M(g,h)$ is the cylinder over
 $[a,b]$. Let $g_i, h_i : X \lo Y, 1 \leq i \leq k$ be maps.
 By the {\em telescope}
 $M((g_1,h_1), \dots, (g_k,h_k))$
 we mean
 the union of the cylinders $M(g_i,h_i)$ where the $h_i$-part of $M(g_i,h_i)$
 is identified with the $g_{i+1}$-part of $M(g_{i+1}, h_{i+1})$ for every $i$.
 Partition $[0,1]$ into $0=t_0< t_1< \dots <t_k=1$
 and consider $M(g_i,h_i)$ as the cylinder  over $[t_{i-1}, t_i]$.
 Then the interval projections of the cylinders define  the interval projection $p$
 from the telescope $M((g_1,h_1), \dots, (g_k,h_k))$
  to the interval $[0,1]$ such that $M(g_i,h_i)=p^{-1}([t_{i-1},t_i])$.
  For pointed maps
 the interval $[0,1]$ is identified with the interval
 in the telescope  so that $[t_{i-1}, t_i]$ is the axis
 of the cylinder $M(g_i,f_i)$. After  such an identification
 we refer to   $[0,1]$ as the telescope axis.
 To shorten the notation  for  $M((g_1,h_1), \dots, (g_k,h_k))$  we  denote
   $\Pi =\{ (g_1,h_1), \dots, (g_k,h_k) \}$
  and then  denote the telescope by
 $M(\Pi)$.  A space is said to be a {\em PL-complex}
 if it admits a triangulation which determines its PL-structure. If
 $X$ and $Y$ are PL-complexes and  the maps in $\Pi$ are PL-maps
  then we always consider $M(\Pi)$ with an induced  PL-structure
  for which each cylinder projection
  $P_i : X \times [t_{i-1},t_i] \lo M(g_i,h_i)$ and the interval projection $p : M(\Pi)\lo [0,1]$
  are PL-maps.
 \\
 \\
 Let $B$ be  an  $l$-dimensional  ball with center $O$. We consider
 $B$ as a pointed space with $O$ as the base point. By a {\em simple map}
 $f : B \lo B$ we mean a map which radially extends a finite-to-one PL-map
 from $ \partial B$  to  $\partial B$.
  More precisely, we consider
 $B$ as the cone over $\partial B$ with the vertex at the center $O$.
Then a simple map is the cone of a  map from $\partial B$ to $\partial B$.
For a simple map $f : B \lo B$ we denote by $\partial f$ the map
$\partial f : \partial B \lo \partial B$ which is the
 the restriction of $f$ to $\partial B$ and by
 the degree of a simple map $f : B \lo B$  we mean  the degree of $\partial f$.
 We say that a simple  map is {\em  non-degenerate } if it has a non-zero degree.
  Clearly  simple maps preserve the base point of $B$.
  Let $g,h: B \lo B$ be simple maps. We consider
  the  cylinder $M(\partial g, \partial h)$ as embedded in $M(g,h)$
  and we will call  $M(\partial g, \partial h)$ the boundary cylinder of $M(g,h)$.
  For  a collection   $\Pi= \{ (g_1,h_1), \dots, (g_k,h_k) \}$
of pairs
of simple maps from $B$  to $B$
we denote
$\partial \Pi=\{ (\partial g_1, \partial h_1) \dots, (\partial g_k, \partial h_k) \}$
 and
 we call    $M(\partial \Pi) \subset M(\Pi)$
  the {\em boundary telescope }  of $M(\Pi)$.
   \\
 \\
   Let $f : X \lo Y$ be a PL-map of PL-complexes. 
   A point $x   \in X$ is said to be a   {\em simple point of degree $d$}
  of the map $f$ if there are closed PL-neighborhoods $X'$ and $Y'$ of $x$ and
  $y=f(x)$ respectively such that such that $f(X')=Y'$, $X'=f^{-1}(Y')$ and
  there are PL-homeomorphisms from
  $X'$ and $Y'$  to a ball $B$
  sending  $x$ and  $y$  to the
  center $O$ of $B$   such that   $f$ restricted to  $X'$ and $Y'$  acts as a simple  map
  of $B$ of degree $d$.  We say that a simple point
  is {\em non-degenerate} if it has non-zero degree. Note that if $f : S^n \lo S^n$
  has a simple point of degree $d$ then $\deg f=d$. Also note that if we have
  a collection of disjoin $n$-balls $B_1, \dots, B_m\subset S^n$ which are $PL$-embedded  in
  $S^n, n \geq 3,$
  and a collection of simple maps $f_i : B_i \lo B_i$ such that $\deg f_i=d$ for every $i$
   then there is a finite-to-one PL-map $f : S^n \lo S^n$ so that $f$ extends $f_i$
  and  $f^{-1}(B_i)=B_i$ for every $i$ and, as  noted before, $\deg f=d$.

 Consider   a telescope
    $ M(\Pi)$    of
   a   collection  $\Pi=\{ (g_1,h_1), \dots, (g_k,h_k)\}$ and
   of non-degenerate simple maps of an $l$-ball $B$, $l \geq 3$ and
  suppose  that $[0,1]$ is partitioned into
   $0=t_0<t_1 <\dots <t_k=1$ such that
     $M(g_i,h_i)=p^{-1}([t_{i-1}, t_i])$ where $p : M(\Pi) \lo [0,1]$ is the interval projection.
  \begin{proposition}
  \label{first-approximation}
    For every map $\phi : X \lo M(\Pi)$ from a space
   $X$ with $\dim_\q X \leq l $
   there is  a  finite-to-one PL-map
 $\psi :   M( \partial \Pi)
 \lo  \partial B \times [0,1] $ such
 that
\begin{itemize}

 \item  $\psi^{-1}(\partial B \times [t_{i-1},t_i])=M(\partial g_i, \partial h_i)$,
  $1\leq i\leq k$;

 \item for every $1\leq i\leq k$ the map $\psi$ has
 non-degenerate simple points   $z_i' ,z_i''\in M(\partial g_i,\partial h_i)$
 such that $z_i' \neq z_i''$, $ p(z'),p(z'') \in (t_{i-1}, t_i)$
 and $\psi(z_i'),\psi(z_i'') \in \partial B \times (t_{i-1},t_i)$;

 \item $\phi$ restricted to $\phi^{-1}( M(\partial \Pi))$
 and followed by $\psi$ extends over $X$ to
 a map $\phi' : X \lo \partial B\times [0,1]$
  such that $\phi'^{-1}(\partial B \times ([t_{i-1},t_i])=\phi^{-1}(M(g_i,h_i))$,
 $1\leq i\leq k$.
\end{itemize}
 \end{proposition}
 {\bf Proof.}  Take simple PL-maps
  $f_i : B \lo B, 1\leq i \leq k+1$
   such that $\deg f_{i} +\deg g_{i}=\deg f_{i+1} +\deg h_{i}$
   and assuming  that $\deg f_1$ is  sufficiently large  we may also  assume
   that $\deg f_{i} +\deg g_{i}=\deg f_{i+1} +\deg h_{i}\neq 0$.
   Denote $g'_i=f_i \circ g_i $, $h'_i=f_{i+1}\circ h_i, 1\leq i\leq k$, and
  $\Pi' =\{(g'_1,h'_1), \dots (g'_k, h'_k) \}$.
 Then there is a natural  projection  $\pi : M(\Pi) \lo M(\Pi')$ and it is
 easy to see that we can replace $\Pi$ by $\Pi'$ and $\phi$ by $ \pi \circ \phi$
 and assume that   $\deg g_i =\deg h_i$ for every  $1\leq i \leq k$.
 \\
 \\
Fix points $z'_i \neq z_i'', z'_i , z''_i \in \partial B \times (t_{i-1}, t_i) $.
 Consider
a finite-to-one PL-map $\alpha_i : \partial (B \times [t_{i-1}, t_i]) \lo  \partial (B \times [t_{i-1}, t_i] )$
so that $\alpha_i$ coincides with $g_i$ on $B \times \{t_{i-1}\}$,
with $h_i$ on $B \times \{t_{i-1}\}$,
$\alpha_i^{-1}(B \times \{t_{i-1}\})=B \times \{t_{i-1}\}$,
$\alpha_i^{-1}(B \times \{t_{i}\})=B \times \{t_{i}\}$,
  $z'_i$ and $z''_i$ are
simple points of $\alpha_i$ of $\deg=\deg g_i =\deg h_i$
and $\alpha_i(z'_i)=z'_i$ and $\alpha_i(z''_i)=z''_i$. Such a map $\alpha_i$ exists because
 $\partial (B \times [t_{i-1}, t_i])$ is a sphere of $\dim \geq 3$.
 Extend $\alpha_i$ to a finite-to-one PL-map
 $\alpha'_i :  B \times [t_{i-1}, t_i] \lo   B \times [t_{i-1}, t_i]$ so that
 ${\alpha'}_i^{-1}(\partial ( B \times [t_{i-1}, t_i])) \subset  \partial (B \times [t_{i-1}, t_i])$.
 Then the cylinder projection $P_i :   B \times [t_{i-1}, t_i] \lo M(g_i,h_i)$
 factors through $\alpha'_i $ and a  finite-to-one PL-map
 $\alpha''_i : M(g_i, h_i) \lo  B \times [t_{i-1}, t_i]$ so that
 $\alpha''_i$ sends the $g_i$-part and the  $h_i$-part of $M(g_i,h_i)$
 by the identity maps to $B \times \{t_{i-1}\}$  and  $B \times \{t_{i}\}$ respectively.
 The maps $\alpha''_i$ define a map $\beta : M(\Pi) \lo B \times [0,1]$.
 We will also  denote by   $z'_i$ and $z''_i$  the images   $P_i (z_i')$  and $P_i(z_i'')$
 of $z'_i$ and $z''_i$ respectively
 in  $M(g_i,h_i)$ under the cylinder projection. Thus we have
 that the points $z'_i$ and $z_i''$ are non-degenerate simple points of $\beta$ restricted to
 $M(\partial \Pi)$ and $\partial B \times [0,1]$.
  Denote  $\phi_0 =\beta \circ \phi : X \lo B \times [0,1]$.
 \\
 \\
 We will construct  for every $ 0\leq j \leq k$ a map $\phi_j : X \lo B \times [0,1]$ and
 a  finite-to-one PL-map $\beta_j : B\times [0,1] \lo B\times [0,1]$  such that for $1\leq i \leq k$
 we have that $\phi_{j+1}$ and $\beta_{j+1} \circ\phi_j$ coincide on
  $\phi_j^{-1}(\partial B \times [0,1])$,
 $\phi_{j+1}^{-1}(B\times [t_{i-1},t_i])=\phi_j^{-1}(B \times [t_{i-1},t_i])$,
   $ 0\leq j \leq k-1$,
 $\phi_j(X)\subset  \partial B \times  [t_0,t_j] \cup B \times [t_j , t_k]$,
 $\beta^{-1}_{j}(B \times [t_{i-1},t_i])\subset B \times [t_{i-1},t_i]$,
 $\beta^{-1}_{j}(\partial B \times [t_{i-1},t_i])\subset \partial B \times [t_{i-1},t_i]$,
 $\beta^{-1}_{j}(B\times \{0\})\subset B\times \{0\}$,
  $\beta^{-1}_{j}(B\times \{1\})\subset B\times \{1\}$, $z'_i$ and $z''_i$ are   non-degenerate
  simple points of $\beta_j$ restricted to  $\partial B \times [0,1]$,
   and $z'_i$ and $z''_i$ are also fixed points of $\beta_j$.
  Set $\beta_0=id$, assume that the construction is completed for $j$ and proceed to $j+1$ as follows.
  \\
  \\
  Denote $X'=\phi_j^{-1} ( B \times \{ t_j \}) $ and $X'' =\phi_j^{-1} ( \partial B \times \{ t_j \}) $.
 By Proposition \ref{rational-dimension} there is a map $g : B \lo B$ of non-zero-degree
 such that if we consider $g$ as a map of  $ B \times \{ t_j \}$ then   the map
 $\phi_j$ restricted to  $X''$ and followed by $g$ continuously extends
 to $\Phi' : X' \lo \partial B \times \{ t_j \}$.
   Once again for every $1\leq i \leq k$
 take a finite-to-one PL-map $\alpha_i : \partial (B \times [t_{i-1}, t_i]) \lo  \partial (B \times [t_{i-1}, t_i] )$
so that $\alpha_i$ coincides with $g$ on $B \times \{t_{i-1}\}$ and
 on $B \times \{t_{i-1}\}$,
$\alpha_i^{-1}(B \times \{t_{i-1}\})=B \times \{t_{i-1}\}$,
$\alpha_i^{-1}(B \times \{t_{i}\})=B \times \{t_{i}\}$
 and $\alpha_i(z'_i) =z'_i$, $\alpha_i(z''_i)=z''_i$ and $z'_i$ and $z''_i$ are
simple points of $\alpha_i$ of  $\deg g$.
Extend $\alpha_i$ to a finite-to-one PL-map
 $\alpha'_i :  B \times [t_{i-1}, t_i] \lo   B \times [t_{i-1}, t_i]$ so that
  ${\alpha'}_i^{-1}(\partial ( B \times [t_{i-1}, t_i])) \subset  \partial (B \times [t_{i-1}, t_i])$.
 and denote by $\beta_{j+1}: B\times [0,1] \lo B \times [0,1]$ the map
 defined by $\alpha'_i$, $1 \leq i \leq k$.
 Let $\Phi''=\beta_{j+1} \circ \phi_j : X \lo B\times [0,1]$.
  Recall that  $\Phi''$ restricted to $X''$  extends
to a map  $\Phi' : X' \lo \partial B \times \{t_j\}$.
Since
$Y=\partial  B \times [t_j, t_{j+1}]\cup B \times \{ t_{j+1} \}$ is contractible
we can extend  $\Phi''$ restricted to $\Phi''^{-1}(Y))$
to a map $\Phi''' : \Phi''^{-1}(B \times [t_j, t_{j+1}]) \lo Y $ so that
$\Phi'''$ coincides with $\Phi'$ on $X'$ and with $\Phi''$ on
$\Phi''^{-1}(Y)$. Thus  we define a map $\phi_{j+1} : X \lo B \times [0,1]$
by changing  $\Phi''$ on ${\Phi''}^{-1}(B\times [t_j,t_{j+1}])$ according to $\Phi'''$ and get  that
$\phi_{j+1}(X) \subset \partial B \times [t_0,t_{j+1}] \cup B \times [t_{j+1} ,t_k]$.
It is easy to see that all the properties of $\phi_{j+1}$ and $\beta_{j+1}$ that we required
are satisfied.
\\
\\
Thus we finally  construct a map $\phi_k$ with
$\phi_k (X) \subset \partial B \times [0,1] \cup B \times \{t_k\}$.
It is obvious that the construction of $\phi_{j+1}$ and $\beta_{j+1}$
also applies to construct from $\phi_k$  a map $\phi_{k+1} : X \lo B \times [0,1]$
and a map $\beta_{k+1}$ so that all the properties of $\phi_j$ and $\beta_j$
will be satisfied except the one regarding $\phi_j(X)$ which we change
to $\phi_{k+1} (X)\subset \partial B\times [0,1]$.
Set $\phi'=\phi_{k+1}$, $\psi=\beta_{k+1} \circ \dots\circ  \beta_0\circ \beta|_{M(\partial \Pi)}$.
Clearly  the construction can be carried out so that
$z'_i$ and $z''_i$ will also  be simple points of $\psi$ and the proof is completed.
 $\black$.
 \\
 \\
   Let $Z$  be a CW-complex with $\dim Z =l+1$.
   Assume that $Z$ also has a PL-structure
   which agrees with the CW-structure. By this we mean that for every closed cell $C$ of $Z$
   both $C$ and $\partial C$ are PL-subcomplexes of $Z$.
   A cylinder $M(g,h)$ of non-degenerate simple maps of an $l$-ball $B$
   is said to be {\em properly embedded} in a closed  $(l+1)$-cell  $C$ of $Z$ if
   $M(g,h) \subset C$ is a PL-embedding into $C$,
   $\partial C \cap M(g,h)=$ the union of the $g$-part and
   the $h$-part of $C$ and $M(g,h)\setminus M(\partial g,\partial h)$ is open in $C$.
   A  telescope  $M(\Pi)$ of a collection
     $\Pi=\{ (g_1,h_1),\dots,(g_k,h_k)\}$   of non-degenerate simple maps of an $l$-ball $B$ is said to be
   {\em properly embedded} into $Z$ if there are   closed  $(l+1)$-cells $C_1, \dots, C_k$
   of $Z$  such that  each cylinder $M(g_i,h_i)$ is properly embedded into $C_i$, $1\leq i\leq k$,
   and $M(\Pi)\setminus M(\partial \Pi)$ is
   open in $Z$.
       An embedding $\gamma: [0,1] \lo Z$ is said  to be {\em a proper  path} from $z_0=\gamma(0)$ to
    $z_1 =\gamma(1)$ if
    the path $\gamma$ identified with $[0,1]$ is the axis of a telescope
    properly embedded  in $Z$. We always assume that a proper path $\gamma$ is parametrized
    by the interval projection $p$ of
    a telescope witnessing that $\gamma$ is proper, that is
    $p(\gamma(t))=t$ for every $t \in [0,1]$.
         It is easy to see that for a proper path $\gamma$ we can  choose
         a telescope $M(\Pi)\subset Z$
    witnessing that $\gamma$ is a proper path to be so close
    to the path $\gamma$ that the interval projection $p : M(\Pi) \lo [0,1]$ to the axis
    will
    have the diameter of the fibers  as small as we wish.

 \end{section}

 \begin{section}{Proof of Theorem \ref{general-exotic}}
\label{proof-general}
 Let $X$ be an $n$-dimensional Boltyanskii compactum, $n\geq 6$.
Recall that  $d_X \leq B_n$.
Consider the dimension types $D_1$ and $D_2$ defined for every prime $p$ by: 
$$ D_1(p)=3^-,\ \ D_1(\q)=2,\ \ \ \ \ \  and\ \ \ \ \ \ D_2(p)=(n-5)^+,\ \ \ D_2(\q)=n-4.$$
Note that $B_n= D_1 \oplus D_2 +1$ and $\dim (D_1+1)\boxplus D_2 =n-1$.
By Theorem \ref{decomposition} there is a decomposition
$X =A \cup B$ of $X$ such that $d_A \leq D_1$ and $d_B \leq D_2$. We may assume
that $A =X \setminus B$ and, by  Theorem \ref{completion}, we may also assume that
$B$ is $G_\delta$ and $A$ is $\sigma$-compact. Represent $A=\cup A_i$ as
a countable union of compact sets $A_i$ such that $A_i \subset A_{i+1}$.

We will construct for each $i$
a  $4$-dimensional compact PL-complex $Y_i$,
a bonding map $\omega^{i+1}_i : Y_{i+1} \lo Y_i$ and
a map $\phi_i : X \lo Y_i$.
 We fix  metrics in $X$ and   in  each  $Y_i$
and with respect to these metrics we  determine  $0<\epsilon_i< 1/2^i$
 such  that
 the following properties will be satisfied:
 \\

 (i) for every open  set $U \subset Y_i$ with $\diam U <2\epsilon_i$
 the set  $\phi_i^{-1}(U)\cap A_i $ splits into at most four  disjoint sets open in $A_i$ and
 of $\diam \leq 1/i$;

\

 (ii) $\dist(\omega^{i+1}_j\circ \phi_{i+1},\omega^{i}_j \circ \phi_{i}) < \epsilon_j/2^{i}$
 for $i>j$ and
  $\dist(\omega^{i+1}_i\circ \phi_{i+1}|_{A_i},\omega^{i}_i \circ \phi_{i}|_{A_i}) < \epsilon_i/2^{i}$
 where $\omega^j_i=\omega^j_{j-1} \circ \dots \circ \omega^{i+1}_i : Y_j \lo Y_i$
for $j>i$ and $\omega^i_i=id : Y_i \lo Y_i$.
 \\
 \\
 The construction will be carried out so that for $Y =\invlim (Y_i,\omega^{i+1}_i)$
  we  have  $\dim_\q Y \leq 3$.
 Let us first show that the theorem follows from this construction.
 Denote  $f_i=\lim_{j\rightarrow \infty}  \omega^j_i\circ\phi_j : X \lo Y_i$.
 From (ii) it follows  that $f_i$ is well-defined, continuous and
 $\dist(f_i, \phi_i) \leq \epsilon_i$ on $A_i$.
 From the definition of $f_i$ it follows that $f_j \circ f^j_i =f_i$.
 Hence the maps $f_i$ define the corresponding  map $f : X \lo Y$ such that
 $\omega_i \circ f=f_i$  where $\omega_i : Y \lo Y_i$ is the projection.
 Then it follows from (i) that for every $y \in Y_i$ the set
 $f_i^{-1}(y)\cap A_i$ splits into  at most four disjoint
 sets closed  in $A_i$ and of $\diam \leq 1/i$. This implies
 that  for every $y \in Y$  we have that $f^{-1}(y) \cap A_i$
 contains at most $4$ points
 and hence $f^{-1}(y) \cap A$ contains at most $4$-points and
 therefore
 $d_f \leq D_2$. Since $\dim Y_i \leq 4$ we have $\dim Y \leq 4$ and
 since $\dim D_2 =n-4$,  Hurewicz Theorem implies that $\dim Y = 4$.
  The condition  $\dim_\q Y \leq 3$ implies that $d_Y(p)\le 4^-$ for all $p$. Therefore,
   $d_Y \leq D_1+1$
 and the theorem follows.
 \\\\
 Let us begin  the construction of $Y_i$, $\omega^{i+1}$ and $\phi_i$.
 In addition to (i) and (ii) we need a few more conditions to be satisfied.
 First we require that\\

 (iii) $\phi_i$ is dimensionally deficient on $A$
 for every $i$. By this we mean that there is triangulation of $Y_i$
 which agrees with the PL-structure of $Y_i$ for which $\phi_i$
 is dimensionally deficient on $A$.
 \\
 \\
 We will endow each $Y_i$ with a structure of a CW-complex which
 agrees with the PL-structure of $Y_i$.
 For every $i$ we  will also construct  a finite collection ${\Gamma}_i$ of
 disjoint proper paths in $Y_i$ and a finite closed cover ${\cal F}_i$ of $Y_i$
  having the following
 properties:
 \\

 (iv) for each $i$ and   each closed $4$-cell $C$  in $Y_i$
     there is a path $\gamma \in {\Gamma}_{i}$ such that
  $\gamma$ crosses $ C$
  (we say that
 a path crosses a closed cell $C$ if it passes through $C \setminus \partial C$);
 \\

 (v)
 for each path
 $\gamma\in \Gamma_i$ we have that
 the union of the closed  $4$-cells of $Y_i$ that $\gamma$ crosses is contained in a set
 of ${\cal F}_i$ and
  the diameters of the sets in  $\omega^i_j({\cal F}_i)$ is less than $ \epsilon_j/2^i$ for $i >j$.
 \\
 \\
 Let  $Y_1$  be a $4$-dimensional simplex.
 By Proposition
  \ref{finite-to-one-approximation},
 take a map  $\phi_1 : X \lo Y_i$  such that $\phi_1$ is
   $4$-to-$1$ on $A_1$.
  Set $\epsilon_1$ to satisfy (i).
 Take a triangulation of
  $Y_1$ into small simplexes, set the CW-structure of $Y_1$ to coincide
  with this triangulation
  and let  ${\Gamma}_1$ contain only one proper path $\gamma$
   to a boundary point of $Y_1$ crossing each simplex of $Y_1$.
   A telescope
  witnessing that $\gamma$ is an proper path is very simple, it is a telescope
   of the identity maps of $B^3$. We assume that the simplexes of $Y_1$ are so small
   that (iv) is satisfied. By Propositions
 \ref{extension} and \ref{finite-to-one-approximation},
 replace   $\phi_1 : X \lo Y_i$ by a map  dimensionally deficient on $A$ and
   so close to $\phi_1$ that (i) remains true. Set ${\cal F}_1 =\{ Y_1 \}$.
   Assume that the construction is completed
  for $i$ and  proceed to $i+1$ as follows.
  \\
  \\
   Fix a number  $\epsilon >0$ which will be
  determined later.  We are going to  refine   the CW-structure
  of $Y_i$  and modify the paths in $\Gamma_i$.
    Take a closed  $4$-cell $C$ of $Y_i$ and  a path $\gamma\in \Gamma_i$
  that $\gamma$ crosses  $C$.
  Let   $M(\Pi)$ be a telescope witnessing that $\gamma$ is a proper path
  and let  $\gamma$ cross
  $\partial C$ at the points $y_\gamma'$ and $y_\gamma''$ with $t'=p(y'_\gamma)$,
  $t''=p(y''_\gamma)$ and $t'< t''$.
   Consider  disjoint small closed neighborhoods of  $y_\gamma'$ and $y_\gamma''$   of the form
 $C_\gamma'=p^{-1}([t', t'+\delta])$  and $C_\gamma''=p^{-1}([t''-\delta , t''])$
 respectively. Then  the neighborhoods $C_\gamma'$ and $C_\gamma''$ are
 PL-subcomplexes of $Y_i$ and taking  $M(\Pi)$  sufficiently close to $\gamma$
 we may assume that  $C_\gamma'$ and $C_\gamma''$ are disjoint
 for distinct  paths in $\Gamma_i$.  Take
 a triangulation of $Y_i$    underlying
 $C_\gamma'$ and $C_\gamma''$ for all  paths $\gamma \in \Gamma_i$
 crossing $C$. Define the new closed $4$-cells
  covering $C$ to be $C_\gamma'$ and $C_\gamma''$ for all  paths $\gamma \in \Gamma_i$
 crossing $C$ and $4$-simplexes not lying in the cells
 $C_\gamma'$ and $C_\gamma''$ we already defined. The cells of $\dim\leq 3$ will be the simplexes of
 $Y_i$ of $\dim \leq 3$. It is easy to see that each path $\gamma$ in $\Gamma_i$
 crossing $C $ can be modified on $C$ so that it will  cross every new closed
 $4$-cell   lying in $C$,  will be proper
 with respect to the new cells  and the modified paths of $\Gamma_i$ will be disjoint.
  This can be done independently
 for every  old closed $4$-cell $C$ of $Y_i$.
  Thus replacing the old CW-structure of $Y_i$ and
 the old collection $\Gamma_i$ by the new ones
 we may assume that every cell of $Y_i$ is of $\diam <\epsilon$
 and (iv) and (v)  remain true  with the cover ${\cal F}_i$ left unchanged.
  \\
 \\
    Take a path $\gamma \in \Gamma_i$ and let
  $M(\Pi)$ be a telescope  of non-degenerate simple  maps of a $3$-ball  $B$ witnessing that
  $\gamma$ is a proper path.
   Then   the interval projection  to the axis
  of $M(\Pi)$ has
  fibers of $\diam < \epsilon$.
  Denote $X'=A_i \cap \phi_i^{-1}(M(\Pi))$
  and $\phi=\phi_i|_{X'}$.
  By Proposition \ref{first-approximation}
  there are a finite-to-one PL-map
  $\psi : M(\partial \Pi) \lo \partial B\times [0,1] $
 and  a map $\phi': X' \lo \partial B\times [0,1]$
  satisfying the conclusions of \ref{first-approximation} for $X$ being replaced by $X'$.
   Define $Y_{i+1}$ to be the quotient
  space of the disjoint union of
  $Y'_i=(Y_i \setminus  M(\Pi))\cup    M(\partial \Pi)$ and $B\times [0,1]$
  obtained by identifying the points of $ M(\partial \Pi) $
  with  $\partial B\times [0,1] $
  by the map $\psi$ and let $q_i : Y'_i \cup B\times [0,1] \lo Y_{i+1}$ be
  the quotient map.
  Note that $B\times [0,1]$ is PL-homeomorphic to
  $(\partial B \times [0,1]\cup B \times \{0\})\times [0,1]$ and consider
   the retraction $ B\times [0,1] \lo \partial B \times [0,1]\cup B \times \{0\}$
  coming from the projection to the first factor in the second
   representation of $B\times [0,1]$.  This retraction defines the corresponding retraction
   $r_{i+1} : Y_{i+1} \lo Y_{i+1}^r=q_i(Y'_i) \cup B\times \{0\}\subset Y_{i+1}$.
    \\
  \\
  Define a map $\phi_{i+1} : X \lo Y^r_{i+1}$ on
  $X\setminus \phi^{-1}_i(M(\Pi) \setminus M(\partial \Pi))$
  as the map $\phi_i$ followed by $q_i$,  on
  $X'$ by $\psi$
  and on $\phi^{-1}_i(M(\Pi) \setminus M(\partial \Pi))$  as any map
  to
  $\partial B \times [0,1]\cup B \times \{0\}$ extending what we already defined.
    Since $\partial B \times [0,1]\cup B \times \{0\}$
  is homeomorphic to $B$  such an extension exists.
  \\
  \\
  Note that by  Proposition \ref{first-approximation} we have that
  the sets  $\partial B\times \{t \}, 0\leq t\leq 1$,
 each fiber of $q_i$  and the sets
 $\phi_{i+1}(\phi^{-1}_i (y)), y \in Y''_i=A_i \cup( Y_i \setminus(M(\Pi) \setminus M(\partial \Pi))$
 are  contained
   in closed $4$-cells of $Y_{i}$. Hence taking $\epsilon$ sufficiently small we can define
   a map $\omega^r_{i+1} : Y^r_{i+1} \lo Y_i$ such that $\omega^r_{i+1} \circ q_i | Y'_i$ will be as close
   to the identity embedding of $Y'_i$ into $Y_i$  as we wish, $\omega^r_{i+1} ( \partial B \times [0,1])$ will be
   as close to the path $\gamma$ as we wish
   and the variation of $\omega^r_{i+1}$ on the sets  $B \times \{0\}$, $\partial B \times \{t\}, 0\leq 1\leq 1$ and
    $\phi_{i}(\phi^{-1}_{i+1} (y)), y \in Y''_i$  will be as small as we wish
   where  the variation  means the supremum of the diameters of the images
   of the sets under $\omega^r_{i+1}$. Set
   $\omega^{i+1}_i =\omega^r_{i+1} \circ r_{i+1} : Y_{i+1} \lo Y_i$. Then we may assume
   that the variation of $\phi_i$ on the fibers of  $\omega^{i+1}_i \circ \phi_{i+1} $ over $Y''_i$
   is as small as we wish and $\phi_{i}(\phi_i^{-1}(M(\Pi))$ is contained in any given
   neighborhood of a set of ${\cal F}_i$ containing all the cells of $Y_i$ that  $\gamma$
   pases through. Thus we may assume that (ii) holds.
   \\
   \\
   Recall that, by Proposition \ref{first-approximation}, $\psi$ is a finite-to-one PL-map.
   Then we define a PL-structure on $Y_{i+1}$ for which $q_i$ is a PL-map.
   Clearly that $q_i$ is finite-to-one and,  since $\phi_i$ is dimensionally deficient on $A$,
   we get, by Proposition \ref{subdivision}, that
   $\phi_{i+1}$ can be replaced by an arbitrarily closed map which is dimensionally
   deficient on $A \cap (Y_i \setminus (\phi_i^{-1}(M(\Pi)\setminus M(\partial \Pi))$.
   Then, by Proposition \ref{extension}, $\phi_{i+1}$ can be again
   arbitrarily closely approximated by a map which coincides with $\phi_{i+1}$ on
   $Y_i \setminus (\phi_i^{-1}(M(\Pi)\setminus M(\partial \Pi))$ and
   dimensionally deficient on $A \cap \phi^{-1}_i(M(\Pi))$. Thus replacing
   $\phi_{i+1}$ by such a map we may assume that (iii) holds.
   By Proposition \ref{finite-to-one-approximation}  we may assume in addition that
   $\phi_{i+1}$ is $4$-to-$1$ on $A_{i+1}$ and determine $\epsilon_{i+1}>0$ for which
   (i) holds.
   \\
   \\
   Since the paths of $\Gamma_i$ are disjoint we can take   telescopes
   witnessing that the paths in $\Gamma_i$ are proper to be disjoint  and
   very close to the paths of $\Gamma_i$. Then
 the above construction can be carried out  simultaneously
 for all the paths
  $\gamma \in \Gamma_i$ and therefore we may assume that constructing $Y_{i+1}$,
 $\phi_{i+1}$ and $\omega^{i+1}_i$ we took care of all the paths of $\Gamma_i$.
 \\
 \\
 Now we will define a CW-structure of $Y^r_{i+1}$, a finite collection
   $\Gamma^r_{i+1}$ of disjoint proper paths in $Y^r_{i+1}$ and establish
   a property that will be used in proving that $\dim_\q Y \leq 3$.
    Let $C$ be  a closed $4$-cell
   of $Y_i$.
    $\gamma \in \Gamma_i$ a path that crosses $C$ and $M(\Pi)$ be the telescope
    witnessing that $\gamma$ is a proper path that we used in the construction.
    Consider the cylinder $M(g,h)$  of the telescope $M(\Pi)$ such that
    $M(g,h)$ lies in $C$. Assume that $M(g,h)$ is a cylinder over the interval $[t',t'']$
    and recall that Proposition \ref{first-approximation} applied $M(\Pi)$ produced
    a map $\psi$  with distinct non-degenerate simple points $z'$ and $z''$
    in $M(\partial g, \partial h)$ such that $ p(z'),p(z'') \in( t', t'')$.
    Note
    that $\psi$ induces the map
   from the sphere $M(\partial g,\partial h)/( M(\partial g,\partial h)\cap\partial C ) $
   to the sphere
   $\partial B \times [t',t'']/(\partial B \times \{t',t''\})$ and, since $z'$ and $z''$ are
    non-degenerate simple points of this induced  map,
    its degree  is equal to the degree of these points
    and hence  different from $0$.
    This implies that\\

    (*) $H^4(q_i(C\cap Y_i'), q_i(\partial C \cap Y'_i);\q)=0$.
    \\
    \\
    Denote $C_r =q_i(C \cap Y_i')$, $z'_r =\psi(z')$ and $z''_r=\psi(z'')$ and
    consider $z'_r$ and $z''_r$ as points of $Y^r_i$.
      Then $z_r'$, $z_r''$ have closed neighborhoods $C'$ and $C'$ in
    $C_r$ of the form
     $C_r'=M(g',id)$ and $C_r''=M(g'',id)$ such that $C_r'\cap C_r''=\emptyset$,
     $C_r'$ and $C_r''$  are disjoint for distinct paths in $\Gamma_i$ that
     cross $C$,
  $g'$ and $g'$ are simple non-degenerate maps of a $3$-ball $B$,
     $M(g',id)$ and $(M(g'',id)$ are PL-embedded in  $C\setminus (M(\Pi)\setminus M(\partial(\Pi))$,
 $z_r'$ and $z_r''$ are   on the axis of $M(g',id)$ and $M(g'',id)$ respectively,
  $M(g',id)\cap \partial B \times (t',t'')$=the $g'$-part of $M(g',id)$ and
   $M(g'',id)\cap \partial B \times (t',t'')$=the $g''$-part of $M(g'',id)$. Take a triangulation
   of $C_r$  which underlies all the sets $C_r'$ and $C_r''$ and
   define  the closed $4$-cells covering $C_r$ to be the sets $C_r'$ and $C_r''$
   and the $4$-simplexes which are not contained in $C_r'$ and $C_r''$
   for all paths of $\Gamma_i$
   that cross $C$. We do that for all closed  $4$-cells $C$ of $Y_i$
   and this way we define the closed $4$-cells of  the  CW-structure
   of $Y_{i}^r$.  As  cells of $\dim <4$ we take the simplexes
   of $\dim <4$ of a triangulation underlying the closed $4$-cells
   we already defined. It is easy to see that for every closed $4$-cell
   $C$ and every path $\gamma \in \Gamma_i$ that crosses $C$ there is a proper path
   $\gamma_r^C$ in $Y_{i+1}^r$ that begins at $z_r'$, ends at $z''_r$ and crosses every
   closed $4$-cell  in $C_r$. Denote by $\Gamma^r_{i+1}$ the collection of the paths
   $\gamma_r^C$ for all closed $4$-cells of $Y_i$. It is clear that taking $\epsilon$
   sufficiently small we may assume that the images of
   $C_r$ under $\omega^r_{i+1}$ are as small as we wish for all closed $4$-cells of $Y_i$.
   Thus we may assume that the images of the cells of $Y^r_{i+1}$ and the paths
   in $\Gamma^r_{i+1}$ under $\omega^r_{i+1}$  are as small as wish.
   \\
   \\
   Now we will extend the CW-structure of $Y^r_{i+1}$ to a CW-structure
   of $Y_{i+1}$, define a finite collection of proper paths $\Gamma_{i+1}$ in
   $Y_{i+1}$ and a finite closed cover ${\cal F}_{i+1}$ of $Y_{i+1}$.
    Take a path $\gamma \in \Gamma_i$
   and let $B\times [0,1] \subset Y_{i+1}$ be the set used in the construction
   for the path $\gamma$. Then every  path $\gamma_r^C$  (for a closed $4$-cell $C$
   of $Y_i$ crossed by $\gamma$) has two end points $z'_r$ and $z''_r$ in $\partial B\times [0,1]$
   and these points are distinct for different cells  $C$ of $Y_i$ crossed by $\gamma$.
   Denote by $Z$ the collection of all  the points $z',z''$ for all the cells $C$ crossed
   by $\gamma$. Take a PL-partition $\cal E$  of $B\times [0,1] \cup \partial B \times [0,1]$
   into closed $3$-cells  such that each point in $Z$ belongs to the interior of
   a cell in $\cal E$. Define the closed $4$-cells covering $B\times [0,1]$ to be
   $r_{i+1}^{-1}(E)$, $E \in\cal E$.  Thus together with the closed $4$-cells of $Y^r_{i+1}$
   we have defined the closed $4$-cells of $Y_{i+1}$ covering $Y_{i+1}$ and as we already
   did before we define cells of $\dim \leq 3$  to be simplexes of $\dim \leq 3$ of
   a triangulation of $Y_{i+1}$ underlying the already defined closed $4$-cells of $Y_{i+1}$.
    Extend every path $\gamma_r^C \in \Gamma^r_{i+1}$ to a path
    $\gamma_{i+1}^C$  by adding to the end points
    $z'_r$ and $z''_r$ of $\gamma_r^C$ two paths going along $r_{i+1}^{-1}(z'_r)$ and
   $r_{i+1}^{-1}(z''_r)$ and connecting $z'_r$ and $z''_r$ with the points
   $r_{i+1}^{-1}(z'_r)\cap B\times \{1\}$ and
   $r_{i+1}^{-1}(z''_r)\cap B\times \{1\}$ respectively. It is clear that $\gamma_{i+1}^C$
   is a proper path of $Y_{i+1}$.  For each  cell $E \in \cal E$ that does not contain
   a point of $Z$ we will take any proper path $\gamma^E_{i+1}$  contained in $r_{i+1}^{-1}(E)$ and connects
   any two interior points of $r_{i+1}^{-1}(E) \cap B\times \{1\}$. Denote by $\Gamma_{i+1}$
   the collection of all the paths $\gamma^E_{i+1}$ and $\gamma_{i+1}^C$ we have constructed
   for
   all paths $\gamma \in \Gamma_i$. It is clear that  every closed $4$-cell of
   $Y_{i+1}$ is crossed
   by a path in $\Gamma_{i+1}$ and assuming that the partition $\cal E$ of
   $B\times [0,1] \cup \partial B \times [0,1]$ is fine enough we can also assume
   that  the images of all the cells in $Y_{i+1}$  and all the paths in $\Gamma_{i+1}$
   under $\omega^{i+1}_i$ are as small as we wish.
   For each path in $\Gamma_{i+1}$ consider the set which is the union
   of the closed $4$-cells intersecting the path  and
   define the cover ${\cal F}_{i+1}$ to be the collection of such sets for
   all the paths in $\Gamma_{i+1}$. Thus we have constructed $\Gamma_{i+1}$
   satisfying (iv) and we may assume that
    ${\cal F}_{i+1}$  satisfies (v). The construction is completed.
   \\
   \\
    The last step of the proof is to verify that $\dim_\q Y\leq 3$.
 Note that for    closed $4$-cells $C$  of $Y_i$,
 the sets  $q_i(C\cap Y'_i)\setminus q_i(\partial C\cap Y'_i)$
are open in $Y_{i+1}^r$,  disjoint for different cells $C$ and
 the  complement in $Y^r_{i+1}$ of the union of these sets is
 a PL-subcomplex of $\dim \leq 3$.  Recall that we may assume that
 the images  of $q_i(C\cap Y'_i)$ under $\omega^r_{i+1}$ are as small as we wish.
 Then we may assume that for the map $r_{i+1}\circ \omega_{i+1} : Y \lo Y^r_{i+1}$
 the fibers of  $r_{i+1}\circ \omega_{i+1}$ and the preimages of the sets
 $q_i(C\cap Y'_i)$ for closed $4$-cells $C$ of $Y_i$ are as small as we wish if
 $i$ is sufficiently large. This  together with (*) implies that
 $H^4(Y,Y';\q)=0$ for every closed subset $Y'$ of $Y$ and hence $\dim_\q Y\leq 3$.
 The proof is completed. $\black$

\end{section}

Alexander Dranishnikov\\
Department of Mathematics\\
University of Florida\\
444 Little Hall\\
Gainesville, FL 32611-8105\\
dranish@math.ufl.edu\\
\\
Michael Levin\\
Department of Mathematics\\
Ben Gurion University of the Negev\\
P.O.B. 653\\
Be'er Sheva 84105, ISRAEL  \\
 mlevine@math.bgu.ac.il\\\\
\end{document}